\documentclass[12pt]{article}
\usepackage{epsfig}
\usepackage{amsmath}
\usepackage{amssymb}
\usepackage{amscd}
\usepackage{booktabs}
\usepackage{pstricks}
\newpsobject{showgrid}{psgrid}{subgriddiv=1,griddots=10,gridlabels=6pt}

\begin{document}
\title
{Unilateral Small Deviations for the Integral \\
of Fractional Brownian Motion }
\date{}
\maketitle
\noindent {\footnotesize G.~MOLCHAN$^{1,2}$, A.~KHOKHLOV$^{1}$} \\
\noindent {\footnotesize \,$^{1}$
{\it International Institute of Earthquake Prediction Theory and Mathematical \\
Geophysics RAS, 79, b2, Varshavskoe sh., 117556, Moscow, Russia \\
E-mail: molchan@mitp.ru, khokhlov@mitp.ru}} \\
\noindent {\footnotesize \,$^{2}$
{\it SAND Group, Galileo bld., The Abdus Salam International Centre \\
for Theoretical Physics, Strada Costiera 11, 34014,
Trieste, Italy }\\
}
\bigskip

\noindent {{\footnotesize \bf Abstract.}
{\footnotesize We consider the paths of a Gaussian random
process $x(t)$, $x(0)=0$ not exceeding a fixed positive level over a large
time interval $(0,T)$, $T\gg 1$. The probability $p(T)$
of such  event is frequently a regularly varying function at $\infty$
with exponent $\theta$.
In applications this parameter can provide information on fractal
properties of processes that are subordinate to $x(\cdot)$. For this
reason the estimation of $\theta$ is an important theoretical problem.
Here, we consider the process $x(t)$ whose derivative is fractional
Brownian motion with self-similarity parameter $0<H<1$.
For this case we produce new computational evidence in favor of the relations
$\log p(T)=-\theta \log T(1+o(1))$ and $\theta =H(1-H)$.
The estimates of $\theta$ are to within 0.01 in
the range $0.1\le H\le 0.9$. An analytical result for the problem in hand is
known for the markovian case alone, i.e., for $H=1/2$. We point
out other statistics of $x(t)$ whose small values have
probabilities of the same order as $p(T)$  in the $\log$ scale.}

\bigskip


\bigskip


\newpage

\section{Introduction}

\bigskip
\noindent The asymptotics of tail probabilities
$P(\max \limits_{\Delta}x(s)>u)$, $u \to \infty$  of a
Gaussian random process $x(t)$ is a problem that has received a
sufficiently complete solution (see, e.g., \cite{LLR},
\cite{Pi}). Considerable progress also occurs in the study
of probabilities of small deviations, that is, of events of the form
$\{ \max \limits_{\Delta}\vert x(s)\vert <\varepsilon \}$ as
$\varepsilon \downarrow 0$ (see, e.g., \cite{Li},  \cite{LS}).
However, there is
practically a complete absence of general approaches to the
analysis of probabilities of unilateral small deviations, i.e., of
probabilities of the form

\begin{eqnarray}\label{qq1}
p(u\vert T)=P(\max \limits_{\Delta \cdot T}(x(s)-x(0))<u)
\end{eqnarray}
for small $u$ or large $T$. Here, $T$ is the similarity ratio and $\Delta$ a
finite closed interval that contains 0. If $x(\cdot )$ is a self-similar
process, that is, finite dimensional distributions of $\{ x(t\Lambda)\}$ and
$\{ \vert \Lambda \vert^hx(t)\}$ are identical for any $\Lambda \ne 0$,
then $p(u\vert T)=p(uT^{-h}\vert 1)$, so that the problems
on the asymptotics of (\ref{qq1}) for the $u$ and $T$ indicated above are
equivalent.

Exponential asymptotics are typical of large and small
deviations of $x(\cdot)$, while power law ones are typical of
unilateral small deviations. Power law asymptotics are rather
popular in physics, since they frequently provide information on
fractal properties of physical processes (see the examples later in
this paper). In this connection one is faced with the task of
calculating the exponent

\begin{eqnarray}
\theta = -\lim\log p(1\vert T)/\log T, \quad T \to \infty.
\label{qq2}
\end{eqnarray}

We shall refer to the exponent $\theta$ for a sequence of events
${\cal A}_T$ in what follows, when
$\log P({\cal A}_T)=-\theta \log T(1+o(1))$ as $T\to \infty$.

There are few explicit estimates of (\ref{qq2}). Molchan (\cite{Mo})
has found the exponent $\theta$ for fractional Brownian motion
(FBM), $b_H(s)$. That one-parameter family of processes is
specified by the requirement $b_H(0)=0$ and by the structure function

\begin{eqnarray*}
E\vert b_H(t)-b_H(s)\vert^2 = \sigma \vert t-s\vert^{2H}, \quad 0<H<1,
\end{eqnarray*}
where $H$ is the self-similarity parameter of the process. When
$H=1/2$, $b_H(s)$ becomes Brownian motion.

According to \cite{Mo}, $\theta =1-H$ when $\Delta \cdot T=(0,T)$ and
$\theta =1$ when $\Delta \cdot T=(-T,T)$. The last estimate is of interest in that it is
independent of $H$. The property in question also remains valid
for FBM with multidimensional time, $t\in R^d$. In that case one
has $\theta =d$, if $\Delta \cdot T=\{ t:\vert t\vert <T\}$.

In mathematical physics one is interested in the parameter
$\theta$  for the integral of fractional Brownian motion (IFBM),
that is, for the process
$x(t)=\int^t_0 b_H(s) ds$, $\vert t\vert <\infty$
\cite{Sib}, \cite{VDFN}. By now, the problem of ~$\theta$ ~has been solved for \, $H=1/2$ \,
only, when $(x(t), \, \stackrel{\cdot}{x}(t))$ \, is a Markov process.
Sinai \cite{Sia} showed that $\theta =1/4$ for $\Delta \cdot T=(0,T)$.
Since paths of $x(t)$ in $(-T,0)$ and $(0,T)$ are
independent, it follows that $\theta =1/2$ for $\Delta \cdot T=(-T,T)$.
The general case $0<H<1$ was studied by the present authors 
(\cite{MK}) both analytically and numerically. Our
analysis suggests the following hypothesis:

\begin{eqnarray}\label{qq3}
\theta = \begin{cases}
1 - H,    & \Delta \cdot T = (-T,T) \qquad \qquad \qquad (a) \\
(1 - H)H, & \Delta \cdot T = (0,T)  \qquad \qquad \qquad \quad (b).
\end{cases}
\end{eqnarray}
The hypothesis (\ref{qq3}a) is corroborated by a related result in \cite{Ha}
 and by numerical calculations. As to (\ref{qq3}b), it has been
confirmed numerically only for the interval $0.1\le H\le 0.6$ to a low
accuracy $(\delta =0.03)$. Here we continue the numerical analysis of
our hypothesis (\ref{qq3}b) by making use of analytical
results derived in \cite{MK}, but radically
modify the evaluation strategy. This allows corroboration of
(3b) in the entire range of $H$ to within $\delta <0.01$. One by-product is
to provide support in favor of the following asymptotics:

\begin{eqnarray*}
P(x(s)x(1)>0, \, \, 1\le s\le T) = O((\log T)^{\alpha (H)}T^{-H(1-H)}),
\quad T \to \infty,
\end{eqnarray*}
where $\alpha (H)$ may have the form $H-1/2$.

\section{Examples}
 Consider a few examples where exponents like (2) are
used.

Let $M(t)=\max \limits_{(0,t)}(x(s)-x(0))$ be the record function of
the process $x(\cdot)$. Levy
has shown (see \cite{IM}) that the record function
of Brownian motion is similar to Cantor's staircase. Its points of
growth make a set $S$ of Hausdorff dimension $\mbox {\rm dim}S=1/2$.
The situation for the FBM family is analogous: $\mbox {\rm dim}S=H$.
This result is based on two fact for FBM in $\Delta \cdot T=(0,T)$:
the exponent is $\theta =1-H$, and

\begin{eqnarray}\label{qq4}
\log [P(G(T)<1)/P(M(T)<1)] = o(\log T),
\end{eqnarray}
where $G(T)$ is the position of the maximum $M(T)$ in the
interval $\Delta \cdot T$ (\cite{Mo}).

Denote by $Z_+$ the first zero of
$b_H(t)$ after the time $t=1$; then, similarly to (\ref{qq4}),

\begin{eqnarray*}
\log [P(Z_+>T)/P(M(T)<1)] = o(\log T).
\end{eqnarray*}
Consequently, the log asymptotics of $P(M(T)<1)$ also
determines the log asymptotics of long excursions of FBM. The
interest in such asymptotics is rather broad in physics and
engineering (see  review \cite{Ma}).

There is a similar problem in geology bearing on the
dimension of fractal sets. Molchan and Turcotte  (\cite{MT}) consider a simple one-dimensional
model of sedimentation in shallow seas. 
The model involves two mechanisms that produce the
sediment: tectonic forces which cause relative rise/fall of sea
level on the one hand and erosion on the other. A sea level fall
brings the upper sedimentary layers above the water and causes
a fast (on the geological time scale) erosion of these. The erosion
process causes  time gaps (unconformities) in datings of the sedimentary
layers. A sea level rise provokes sedimentation. The latter
process is considered to be rapid, so that the sediment follows
sea level variations practically continuously. Supposing the sea
level history to be described by a process $x(t)$, the drilling data
(layer depth and date of formation) must be described by the function

\begin{eqnarray*}
y(t) = \min \limits_{s>t}x(s).
\end{eqnarray*}

Points of growth of $y(t)$ (the set $S$) correspond to the dates of
those layers which have been preserved in the sedimentary
sequence during the entire history of sedimentation. The
common assumption makes $x(s)$ a process with stationary
increments; the self-similarity of $S$ in an extensive range of
scales suggests self-similarity for the stochastic component of $x(\cdot)$.
When $x(s)$ is also assumed to be Gaussian, one arrives at the
model $x(s)=as+b_H(s)$, $a>0$ (see e.g. \cite{EM}).
The dimension of the support of $y(t)$ for the above model
of $x(\cdot)$ can be found similarly to the preceding problem, giving
the result $\text {dim}\,S=H$ (\cite{MT}).

The second example is concerned with the fractal nature of
solutions to the Burgers equation with random initial data:

\begin{eqnarray}\label{qq5}
u_t+uu_x&=&\nu u_{xx} \\
\nonumber
u(0,x)&=&v(x), \quad \vert x\vert <\infty
\end{eqnarray}
where $\int^xv(a) \, da=o(x^2)$, $x\gg 1$
(see, e.g., \cite{Wo}).

When the viscosity $\nu$ is infinitely small, equation (\ref{qq5})
when considered in the one-dimensional case describes the
following dynamics of adhesive particles. A particle at the initial
time at position $a$ has the mass $da$ and the momentum
$v(a)da$. The
particle starts moving at velocity $v(a)$ and preserves that
velocity until  the first collision with its neighbors. The colliding
particles stick together and continue their motion by following
the laws of preservation of mass and momentum, and so on. The
question that arises is what is the dimension of the initial
positions of the particles $S$ that have not collided until time
$t_0$. Such positions are called Lagrangian regular points.

The question is amenable to a purely geometrical interpretation. Let
$U(x)$ be a convex minorant of the curve $\xi (x)=\int^x_0v(a) \, da+x^2/2$.
Then the right-hand derivative $U'(x)$ is nondecreasing function, and its points of growth correspond to
Lagrangian regular points  $S$. The dimensionality problem of $S$
in the case $v(a)=b_{1/2}(a)$ has been solved by Sinai 
\cite{Sib}: $\text {dim}\,S=1/2$. The solution is based on estimating
the exponent (2) for
the integral of Brownian motion in the interval $(0,T)$. For the
case $v(a)=b_H(a)$ it can be asserted that $\text {dim}\,S=H$, if
$\theta$  for IFBM in $\Delta \cdot T=(-T,T)$
admits of the bound $\theta \ge (1-H)$ (\cite{MK}).

The exponents (\ref{qq2}) for the maximum of IFBM are closely
related to the exponents of other statistics.
To be more precise we introduce the
following notation: $M(T)=\max \limits_{\Delta \cdot T}x(t)$, $G(T)$
is the position of $M(T)$ in $\Delta \cdot T$,
and $A(T)=\int \limits_{\Delta \cdot T}{\bf 1}_{x(s)>0}\,ds$
is the occupation time of $x(t)$ above 0 in $\Delta \cdot T$.

{\bf Statement 1.} (Molchan \& Khokhlov, \cite{MK}). Let
$\Delta \cdot T=(0,T)$ or $(-T,T)$. For $x(t)=\text {IFBM}(t)$,
$t\in \Delta \cdot T$ the exponents of the events
$\{ M(T)<1\}$, $\{ \vert G(T)\vert <1\}$, $\{ A(T)<1,\, \vert G(T)\vert <T\}$,
$\{ x(t)<0, \,t\in \Delta \cdot T,\, \vert t\vert >1\}$,

(a) are identical, when they exist;

(b) exist or do not exist simultaneously.

It is a known fact that $G(1)$, $A(1)$ and the last zero,  $Z$, for
Brownian motion in (0,1) have identical distributions, namely,
the arcsine law.
Statement 1 is also true for FBM (\cite{Mo}). Therefore
we can consider this result as a weak
version of the arcsine law for FBM and IFBM.
Statement 1 provides a certain
degree of freedom in numerical analyses of (\ref{qq2}).

\section{Evaluation of $\theta$}
 We are going to evaluate the exponent of the event $\{ M(T)<1\}$
for the IFBM process in the interval $(0,T)$, $T\gg 1$ using the Monte
Carlo method. To select a suitable strategy note the following.
Statement 1 gives some
information on the possible structure of a typical IFBM path with a low
maximum, $M(T)<1$. Namely, the position of the maximum $G(T)$ and the
total time $A(T)$ where $\text {IFBM}>0$ do not
practically grow (the growth is most likely to be a logarithmic
one). In most of the cases the path goes to the lower half-plane
after the lapse of a fixed time, because there are no limitations on the
amplitude from below. Consequently, the essential
information concerning the low maximum is available around
the initial point $t=0$.

Suppose we generate IFBM on a uniform lattice
$(0,\delta,\ldots,L\delta)$ using
triangular factorization of the correlation matrix. We note that
the correlation structure of the process should be reproduced
exactly in the case under consideration, since we are dealing
with rare events $\{ M(T)<1\}$. Hence we shall need a memory of order
$L^2$ for generating a large number of samples. However, the
second-degree growth in memory leads to limitations on lattice
width, hence on the information concerning the low maximum.
This heuristic argument is borne out in \cite{MK}. Following the above strategy, we have not
succeeded in evaluating $\theta$ for $0.6\le H\le 1$.

The way out consists in considering the IFBM in a log time scale.
To do this, we consider the Lamperti transformation which converts
$\{ \text {IFBM}(t), t>0\}$ into a stationary process $x(t)$:

\begin{eqnarray}\label{qq6}
x(t) = c \, \exp(-(1+H)t)\text {IFBM}\,(e^t), \quad \vert t\vert <\infty.
\end{eqnarray}

We have normalized the process so as to make $Ex^2(t)=1$, hence
$c=(1+2H)^{-1/2}$. The correlation function $r(t)$ for $x(t)$ is

\begin{eqnarray}
r(t) = 2c^2(1+H)\cosh(Ht)-c^2\cosh((1+H)t)+0.5c^2\vert \sinh (t/2)\vert^{2H+2}
\nonumber
\end{eqnarray}
and

\begin{eqnarray}
r(t) = \begin{cases}
1-0.5(1-H^2)t^2+0.5c^2t^{2+2H}+O(t^4),    & t\to 0  \\
c(H)\cdot \exp(-\rho t)\cdot (1+o(1)), & t\to \infty,
\nonumber
\end{cases}
\end{eqnarray}
where $\rho =\min(H, 1-H)$.

We can see that $x(t)$ has a smoothness of order
$1+H-\varepsilon$ ($\varepsilon >0$), the number of zeroes of
$x(\cdot)$ is locally finite, and the
mean interzero distance is, according to Rice, given by

\begin{eqnarray}\label{qq7}
\Delta_0 = \pi (1-H^2)^{-1/2}.
\end{eqnarray}

The transformation (\ref{qq6}) does not preserve the point ($G(T)$, $M(T)$)
as an extreme one in a sample of $x(\cdot)$. However, the event

\begin{eqnarray}
\{ \text {IFBM}(t)<0,\quad 1<t<T\} = {\cal A}_T
\nonumber
\end{eqnarray}
is easily rewritten in terms of $x(\cdot)$ on a finite interval:

\begin{eqnarray}
{\cal A}_T = \{ x(\tau)<0,\quad 0\le \tau <T'=\ln T\}: = \tilde{{\cal A}}_{T'}.
\nonumber
\end{eqnarray}

Consequently, the desired exponent

\begin{eqnarray}
\theta = -\lim_{T'\to \infty} \ln P(\tilde{{\cal A}}_{T'})/T'
\label{qq8}
\end{eqnarray}
is converted to exponential from a power-law one. As a result,
the evaluation of $\theta$ splits into three steps:

$\bullet$ generating a stationary process $x(t)$, $t>0$;

$\bullet$ evaluating the distribution of the first zero,  $Z$, for $x(t)$;

$\bullet$ finding (\ref{qq8}) for the tail of the distribution of $Z$.

{\it The Generation of} $x(t)$. The process $x(t)$ was generated as a
stationary sequence $x(\delta \cdot k)$ with the exact correlation function
$r(\delta k)$, $k=0,1,...,L$. This was done using the triangular representation

\begin{eqnarray}
x(k\delta ) = \sum^k_{i=0} a(i\vert k)\varepsilon_i
\label{qq9}
\end{eqnarray}
in terms of the standard white noise $\{ \varepsilon_i\}$. The
representation is implemented by using the progressive Schur algorithm (see
\cite{AC}). The discretization step $\delta$ is
specified by the number of points $n_0$ per mean period $\Delta_0$
for the zeroes of $x(t)$ (see (\ref{qq7})). Since $x(\cdot)$ has a smoothness of
order $\sim 1+H$, and we are interested in long excursions of $x(t)$,
we can well use moderate numbers for $n_0$; we had $n_0=50$
in the calculations. The length $L$ can be found from the
requirement $P(Z>L\delta )=\varepsilon$ where $\varepsilon$ is small.

\begin{figure}[h]
\centerline{
\begin{pspicture}(0,0)(10,5)
\psset{dimen=inner}
\rput(5.0,2.){
\begin{minipage}{0.49\linewidth}
\centerline{
\epsfig{file=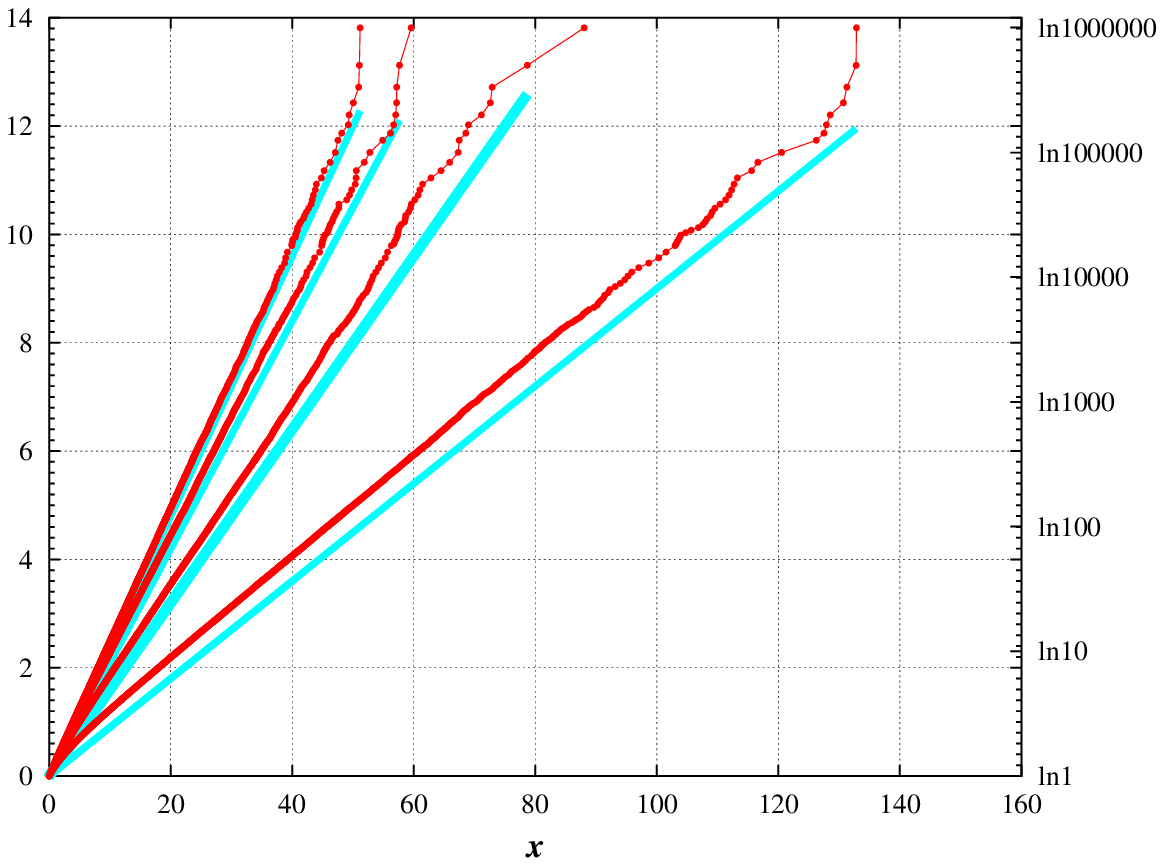, angle=0, width=0.99\linewidth}
}
\end{minipage}
\begin{minipage}{0.49\linewidth}
\centerline{
\epsfig{file=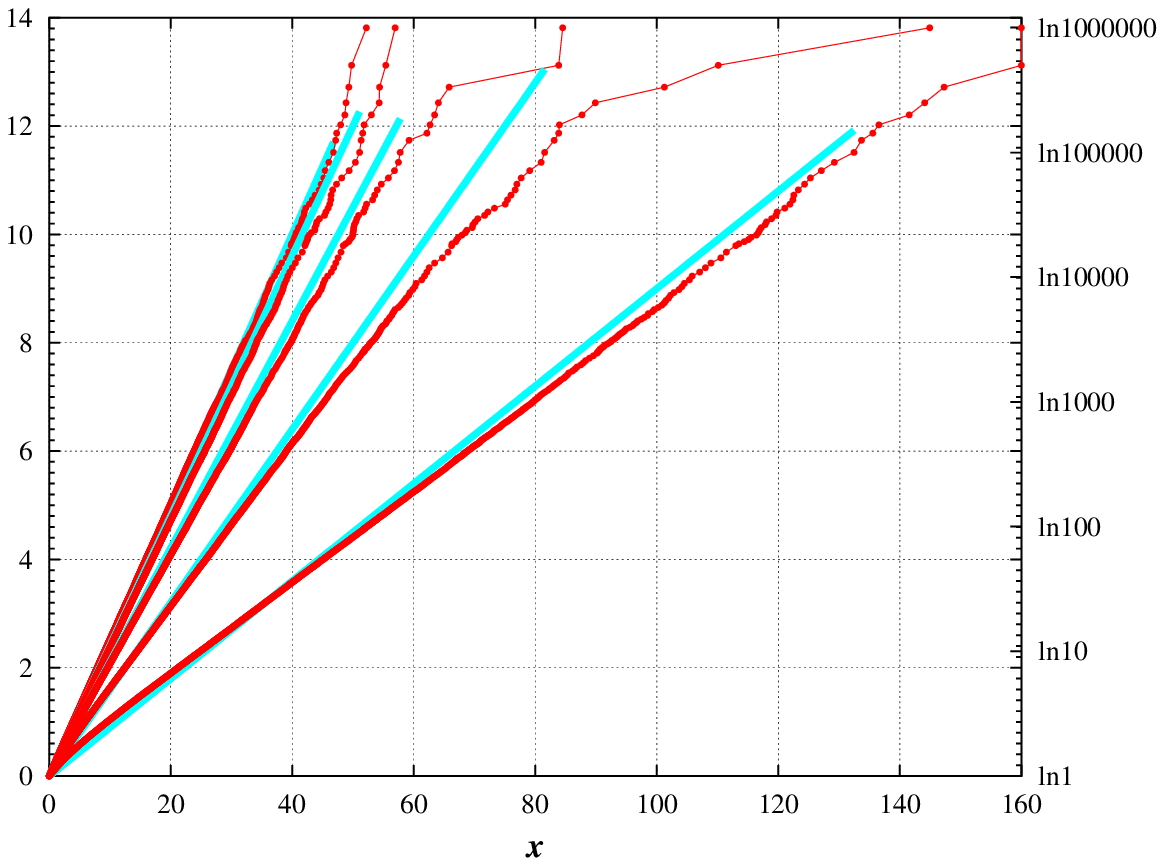, angle=0, width=0.99\linewidth}
}
\end{minipage}
}
\rput(0.2, 3.8){\tiny{0.4}}
\rput(0.6, 3.8){\tiny{0.3}}
\rput(1.1, 3.8){\tiny{0.2}}
\rput(2.9, 3.8){\tiny{0.1}}
\rput(6.9, 3.8){\tiny{0.5}}
\rput(7.3, 3.8){\tiny{0.6}}
\rput(7.7, 3.8){\tiny{0.7}}
\rput(8.5, 3.8){\tiny{0.8}}
\rput(10.2, 3.8){\tiny{0.9}}
\end{pspicture}
}
\caption  {Empirical functions $y=-\ln\,P(Z>x)$ based on
$N=100,000$ paths of IFBM: 
$H=0.1\div 0.4$ (left) and $H=0.5 \div 0.9$ (right).
{\it Straight} lines correspond to  $y=H(1-H)x$.}
\label{pic1}
\end{figure}

As is shown by
some preliminary evaluations of the distribution of $Z$
(see Fig.~\ref{pic1}), the function $P(Z>t)$ for $H=0.1\div 0.9$ is well fitted
with an exponential function: $\exp(-H(1-H)t)$. Consequently, one has
approximately $L\delta \simeq Z_\varepsilon$, where

\begin{eqnarray}
Z_\varepsilon = -\ln \varepsilon /[H(1-H)].
\label{qq10}
\end{eqnarray}
When $\varepsilon =10^{-4}$ and $H=0.2-0.8$, one has $L\delta \simeq 60$,
and $L\delta =100$ when
$H=0.1,\, 0.9$. The Schur algorithm is a recursive one, so it may
become unstable as $L$ increases. The instability manifests itself
in parasitic oscillations of the $a(i\vert k)$ at large $i$
(see (\ref{qq9})). These effects are
typical of $i\cdot \delta \ge 40$ and hardly can  always be overcome by using
available accuracy (for instance {\tt long double} in {\bf C}) . Hence it follows that

$\bullet$ there are computational difficulties in the way of analyzing the
distribution of $Z$ when $H$ is close to either 0 or 1; more exactly,
when $\vert H-0.5\vert >0.4$;

$\bullet$ evaluation of $\theta$ from the values $Z>Z_\varepsilon$,
$\varepsilon <10^{-4}$ calls for higher computation accuracy.

The last conclusion is important, because we do not know
when the log linear asymptotics for the tail of the distribution of
$Z$ becomes valid.

{\it Evaluation of} $\theta$. The parameter $\theta$ was evaluated in the
series of intervals
$\Delta_\varepsilon =(Z_\varepsilon , Z_{\varepsilon /10})$ with
$\varepsilon =0.01; 0.003; 0.001$ for the range $H=0.1-0.9$. Figure~\ref{pic1} does not contradict the assumption of linearity for the plot of
$(t, \ln P(Z>t))$, $t\in \Delta_\varepsilon$. For this reason
we use for the slope of the plot the
maximum likelihood (ML) estimate $\widehat{\theta}$ corresponding to the
distribution $\{ c(\theta) e^{-\theta t}, t\in \Delta_\varepsilon \}$.
Namely, $\widehat{\theta}=x/\vert \Delta_\varepsilon \vert$ where
$\vert \Delta_\varepsilon \vert =\frac{1}{\theta_0}\ln 10$  is
the length of $\Delta_\varepsilon$, $\theta_0=H(1-H)$, and $x$ is the root of

\begin{eqnarray}
x^{-1} - (e^x-1)^{-1} =
[<Z>_\varepsilon - Z_\varepsilon ]/\vert \Delta_\varepsilon \vert.
\label{qq11}
\end{eqnarray}
Here, $<Z>_\varepsilon$ is the empirical mean of all $Z$ observed in
the interval $\Delta_\varepsilon$.

\begin{figure}[h]
\centerline{
\epsfig{file=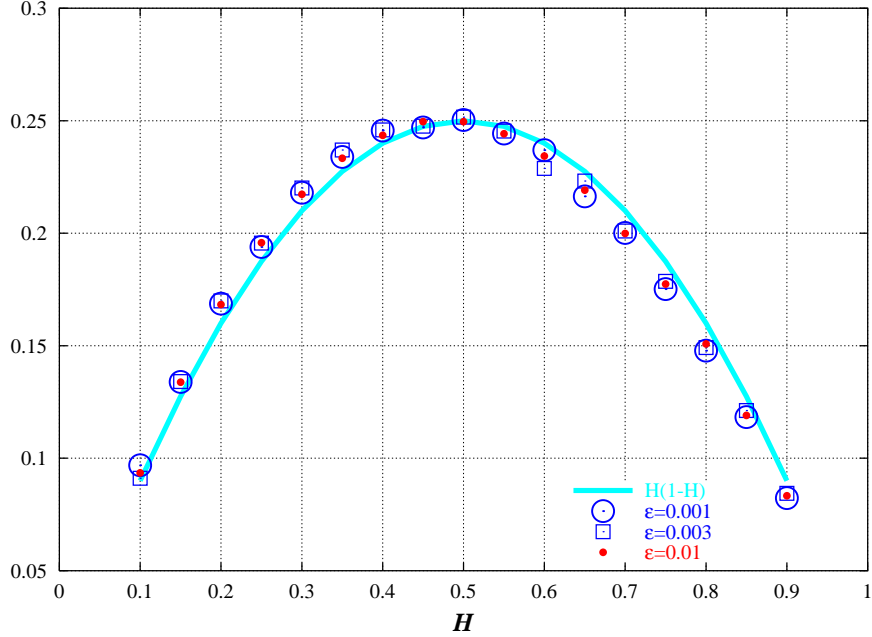, angle=0, width=0.9\linewidth}
}
\caption{Hypothetical $\theta_0=H(1-H)$ and empirical
$\widehat{\theta}$ exponents for three intervals $\Delta_\varepsilon$,  $\varepsilon =0.01,\  0.003,\  0.001$.
Number of paths of IFBM is $N=16\times 300,000$.}
\label{pic2}
\end{figure}

ML estimates of $\theta$ based on
$N=16\times 300,000$ paths of $x(\cdot)$ are shown in Fig.~\ref{pic2}. They
demonstrate that {\it the slopes} $\widehat{\theta}_\varepsilon$ {\it in the intervals}
$\Delta_\varepsilon$, $\varepsilon =0.01, 0.003, 0.001$ {\it are well consistent
among themselves and are identical
with the hypothetical values} $\theta_0=H(1-H)$ {\it to within} 0.01.
Note that the
left-hand endpoint of $\Delta_\varepsilon$ is approximately identical with the
$(1-\varepsilon )$ quantile of the distribution of $Z$. Consequently, the
number of observations $N_\varepsilon$ used to estimate
$\theta_\varepsilon$ is approximately equal to $N\cdot \varepsilon$.

The \, ML \, estimate \, of \, $\theta$ \, for \, the truncated \, exponential
distribution  $\{ c(\theta ) \exp(-\theta t),\,
t\in (\ln \varepsilon^{-1}, \, \ln 10/\varepsilon) \, \theta^{-1}\}$
has the standard deviation
$\sigma_\varepsilon \simeq 1.7014 \, \theta N_\varepsilon^{-1/2}$, where
$N_\varepsilon \gg 1$ is the number of observations. In our case the slope
estimates are close to $\theta_0=H(1-H)$, while $N_\varepsilon$ is large,
hence $\widetilde\sigma_\varepsilon = 1.7 \, \theta_0 N_\varepsilon^{-1/2}$
can serve as a satisfactory theoretical estimate of the standard
deviation for $\widehat{\theta}_\varepsilon$. This statement is
corroborated by our experiments carried out to check the operation
of the random digit generator and the Schur algorithm. The
estimates $\widehat{\theta}_\varepsilon$
were derived above by averaging over 16 serial estimates
$\widehat{\theta}_\varepsilon^{\rm ser}$. Each series consists of 300,000 paths
of $x(\cdot)$. The empirical variance of the averaged estimate
$\widehat{\theta}_\varepsilon$ is in good
agreement with the theoretical value. Examples are given below
in Table~\ref{tab07}.
\begin{table}
\centerline{
\begin{tabular*}{0.6\hsize}{rc@{\quad}|@{\qquad}ccc@{}}
\toprule
 $n_0$ & $\mathbf\varepsilon$ & \Large$\widehat{\theta}$ & \LARGE$\widehat\sigma_{\widehat{\theta}}$ & \LARGE$\widetilde\sigma_{\widehat{\theta}}$\\
\bottomrule
   & .100   &  .2002 & .00009 &  .00011 \\ 
30 & .010  &  .1998 & .00031 &  .00032 \\ 
   & .001 &  .1997 & .00089 &  .00096\\ 
\hline
    & .100   &  .2001 & .00011 &  .00012 \\ 
100 & .010  &  .2001 & .00040 &  .00037 \\ 
    & .001 &  .1989 & .00107 &  .00112\\ 
\bottomrule
\end{tabular*}
}
\caption{ Estimates $\widehat{\theta}$ in intervals $\Delta_\varepsilon$
for $H=0.7$ and 
discretization parameter $n_0 =30$ (top) and  $n_0 =100$ (bottom). 
The standard deviations of $\widehat{\theta}$ are  $\widehat\sigma_{\widehat{\theta}}$ (empirical)
 and $\widetilde\sigma_{\widehat{\theta}}$ (theoretical).
}\label{tab07} 
\end{table}

\begin{figure}[h]
\centerline{
\epsfig{file=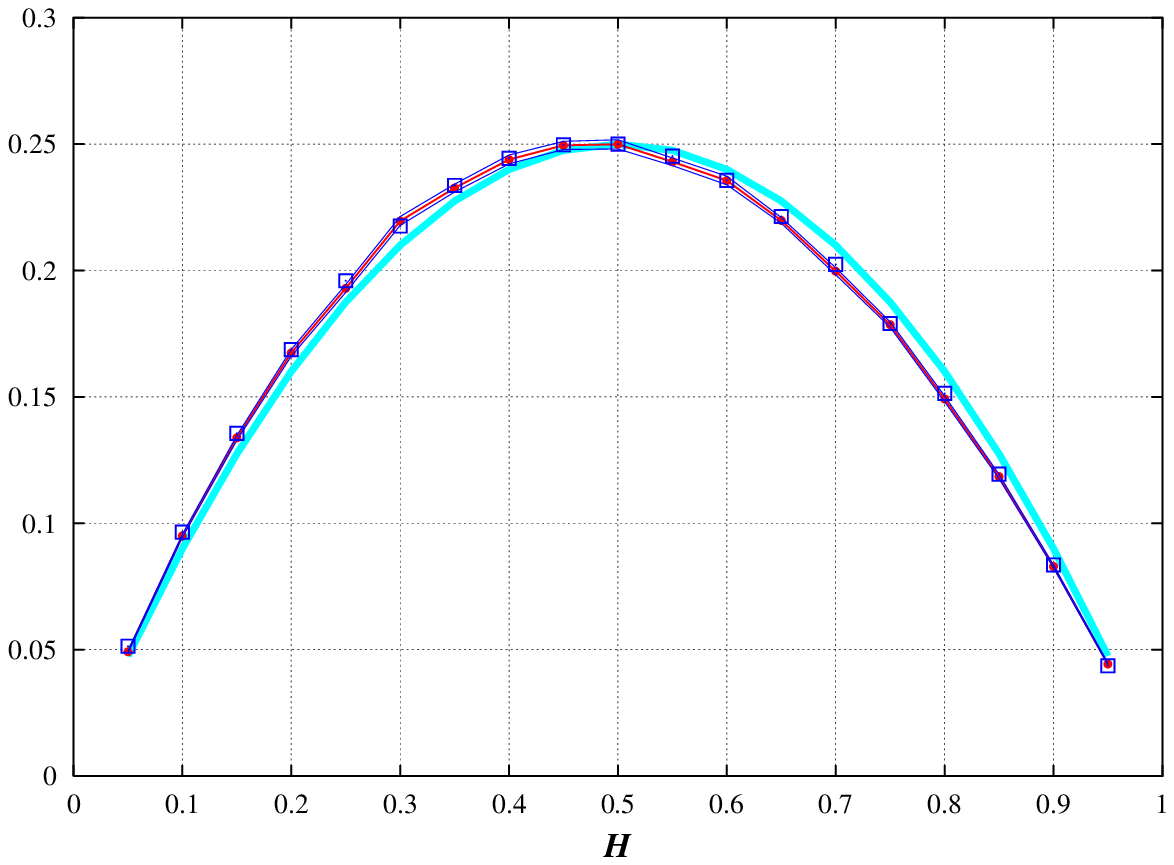, angle=0, width=0.9\linewidth}
}
\caption{Hypothetical values of $\theta_0=H(1-H)$ and interval
estimates $\widehat{\theta}_\varepsilon \pm \widetilde\sigma_\varepsilon$ of the
exponent $\theta$ for the interval $\Delta_\varepsilon$, $\varepsilon=0.01$.
Number of paths is $N=25\times 300,000$.
Boxes correspond to the expected values
$\tilde{\theta}_\varepsilon$ of $\widehat\theta$ under condition \ref{qq12}
and $\alpha =H-0.5$.}
\label{pic3}
\end{figure}

Figure~\ref{pic3} shows interval estimates
$\widehat{\theta}_\varepsilon \pm \widetilde\sigma_\varepsilon$ of the slope
$\widehat{\theta}_\varepsilon$ for $\varepsilon =0.01$. Even though
we have seen above that
$\vert \widehat{\theta}_\varepsilon -\theta_0\vert <0.01$,
Fig.~\ref{pic3} provides evidence of a significant discrepancy
between empirical and hipothetical estimates of $\theta$.
Furthermore we can see that
$\widehat{\theta}_\varepsilon (H)>\theta_0(H)$
when $H<0.5$ and
$\widehat{\theta}_\varepsilon (H)<\theta_0(H)$
when $H>0.5$ for all $\varepsilon =0.01,\,0.003$ and
0.001 (Fig.~\ref{pic2}).

The above inference cannot be ascribed to the effect of
discretization. This is confirmed by the estimates of $\theta$ for
$H=0.7$ with different values of the discretization parameter: $n_0=30$
and 100 steps per period $\Delta_0$ (see Table~\ref{tab07}). The estimates of
$\theta_\varepsilon$ correspond to the intervals
$\Delta_\varepsilon = (\ln \frac{1}{\varepsilon},\,
\ln \frac{10}{\varepsilon})\cdot \theta_0^{-1}$
with $\varepsilon =0.1,\,0.01,\,0.001$
and were derived by averaging over $r$ serial estimates
$\theta_\varepsilon^{\rm ser}$. Each series consists of $10^6$ paths
of $x(\cdot)$, $r=120$ in the case  $n_0=30$ and
$r=90$ in the case  $n_0=100$.
Table~\ref{tab07} also lists empirical and
theoretical standard deviations of the estimates of $\theta_\varepsilon$.
It appears from Table~\ref{tab07} that the estimates
$\widehat\theta_\varepsilon \approx 0.200$ are independent of the
discretization step and there are significant deviation of
$\widehat\theta_\varepsilon$ from the
hypothetical value $\theta_0=0.21$.

We thus have to reject the hypothesis that the tail of the
distribution of $Z$ has a purely exponential asymptotics with
parameter $\theta =\theta_0$. However, this does not rule out the hypothesis
proper of the exponent $\theta =\theta_0$ for the events
$\{ Z>T\} ,T\gg 1$. We are going to show that the distribution

\begin{eqnarray}
P(Z>t) = ct^\alpha e^{-\theta_0t}(1+o(1)), \quad t \to \infty
\label{qq12}
\end{eqnarray}
is consistent with our estimates of $\theta_\varepsilon$.

\begin{figure}[h]
\centerline{
\epsfig{file=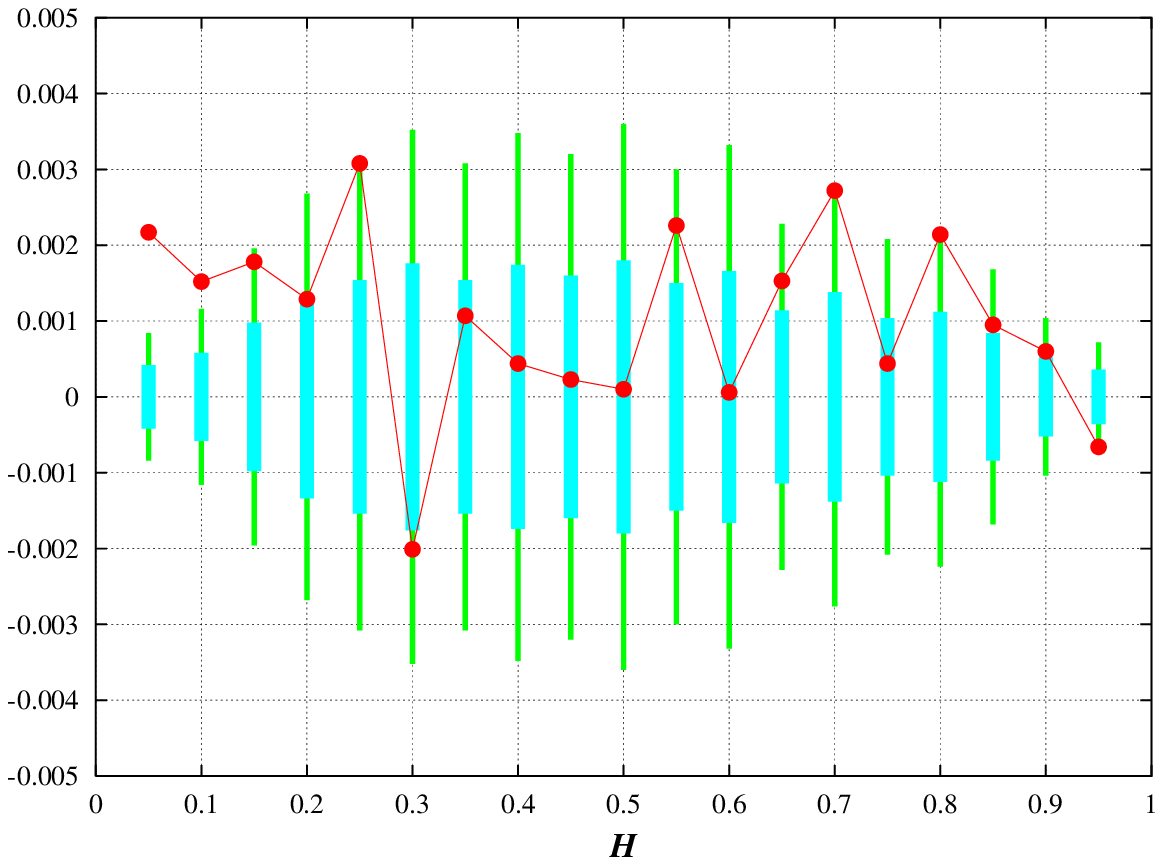, angle=0, width=0.9\linewidth}
}
\caption{Residuals
$R = \tilde{\theta}_\varepsilon -\widehat{\theta}_\varepsilon$,
$\varepsilon =0.01$ (see  Fig.~\ref{pic3}).
{\it Vertical  lines} correspond to two levels of
$\widehat{\theta}_\varepsilon$ deviations:
a)$\pm \widehat\sigma_\varepsilon$ ({\it bold}) and b)$\pm 2\widehat\sigma_\varepsilon$,
({\it thin}) where $\widehat\sigma_\varepsilon$ is empirical standard deviation of
the estimate $\widehat\theta$}
\label{pic4}
\end{figure}

To do this, let us replace the empirical mean $<Z>_\varepsilon$ in (\ref{qq11})
with $EZ{\bf 1}_{z\in \Delta_\varepsilon}$ assuming (\ref{qq12}) and
$o(1)=0$. Solving the equation yields the expected value of
$\widehat{\theta}_\varepsilon$ under (\ref{qq12}).
These estimates $\tilde{\theta}_\varepsilon$  are shown
in Fig.~\ref{pic3} for $\alpha =H-0.5$. It is
seen that the $\tilde{\theta}_\varepsilon$ are in very good agreement
with the empirical estimates  $\widehat{\theta}_\varepsilon$ for
all $H=0.1-0.9$. Figure~\ref{pic4} provides a more detailed view of the residuals
$R_\varepsilon =\tilde{\theta}_\varepsilon -\widehat{\theta}_\varepsilon$.

Both Figs.~3 and 4 show that the empirical estimates
$\widehat{\theta}_\varepsilon$ can
be well fitted using the extra parameter $\alpha$.
However, it is very difficult to get $\alpha$ with  a suitable
resulution $\delta$, say $\delta =0.02-0.03$, if $\alpha$ is small, as
is the case for the model $\alpha =H-0.5$. Indeed, suppose the
parameter $\theta$ is known and $\theta =\theta_0$. The Cramer-Rao
inequality yields the optimal variance $\sigma^2_{\rm opt}$ for $\alpha$:

\begin{eqnarray}
\sigma_{\rm opt} \simeq N_\Delta^{-1/2}\sigma^{-1}(\ln Z_\Delta)
\nonumber
\end{eqnarray}
where $\sigma^2(\ln Z_\Delta)$ is the variance of $\ln Z$ for observations
of $Z$ in the interval $\Delta$, $N_\Delta$ being the number of the observations.
Hence

\begin{eqnarray}
N_\Delta \simeq [\sigma_{\rm opt} \cdot \sigma (\ln Z_\Delta)]^{-2}.
\nonumber
\end{eqnarray}
If $\Delta =(Z_\varepsilon,\,A/\theta_0)$ and $Z_\varepsilon$ is given
by (\ref{qq10}), then the total number of paths of $x(\cdot)$ is

\begin{eqnarray}
N \simeq [\sigma_{\rm opt} \cdot \sigma (\ln Z_\Delta)]^{-2}
\varepsilon^{-1}.
\nonumber
\end{eqnarray}
One has $\vert \alpha \vert <0.5$ in the model $\alpha=H-0.5$. For this
range of $\alpha$ with $\varepsilon =0.01$ and 0.001 for $A=20$ one has
$\sigma (\ln Z_\Delta) \simeq .15883$ and
0.11536, respectively. The requirement $\sigma_{\rm opt}=0.01$ makes the
number $N$ large enough: $40\cdot 10^6$ if $\varepsilon =0.01$ and
$750\cdot 10^6$ if $\varepsilon =0.001$.

\section{Conclusion}

We have shown that the exponent hypothesis  $\theta_0=H(1-H)$
for a series of statistics related to IFBM (see Statement~1) is well
corroborated by our computations. Further refinement of the
tail probabilities for these statistics faces considerable computational
difficulties in view of the amount of computation required, the computation
accuracy, and checks on the random digit generator.

The symmetry of the exponent: $\theta (H)=\theta (1-H)$  which
$\theta_0$ has by definition
is not obvious in our problem because of quite different
properties of FBM for  $H<1/2$ and $H>1/2$. Our analysis suggests
that this difference can manifest itself in more sophisticated asymptotics,
e.g.

\begin{eqnarray}
P({\rm IFBM}(t)<1, \,\, 0<t<T) = O(T^{-\theta_0}(\log \,T)^{\alpha(H)})
\nonumber
\end{eqnarray}
where $\alpha(1/2)=0$. Of course this relation, as well as $\theta_0=H(1-H)$,
are needed in analytical corroboration.

\section*{ Acknowledgments}

This research was supported by the Russian Foundation for
Basic Research (grant 02-01-00158).


\end{document}